\ifcase0
   \documentclass[12pt,a4paper]{article}
   \usepackage{euler,amsmath,amsfonts,url,defns}
   \usepackage[colorlinks]{hyperref}
\or\documentclass[11pt,a5paper]{article}
   \usepackage{euler,amsmath,amsfonts,defns,ajr}
\or\documentclass[12pt,a4paper]{refart}
   \usepackage{euler,amsmath,amsfonts,defns,ajr}
   \settextfraction{0.80}
\fi

\usepackage{ProTex}
\AlProTex{red,<<<>>>,title,list,`}

\allowdisplaybreaks

\newcommand{\ibc}{\textsc{ibc}s}

\newcommand{\uij}{u_{i,j}}

\newcommand{\bdelta}{\text{\boldmath$\delta$}}

\title{Computer algebra derives the slow manifold of patch or element dynamics on lattices in two dimensions}

\author{Tony MacKenzie\thanks{Department of Mathematics and Computing, University of Southern Queensland, Toowoomba, Queensland~4352, Australia.} 
\ and 
A.~J. Roberts\thanks{Corresponding author: School of Mathematical Sciences, University of Adelaide, South Australia~5005, Australia. 
\protect\url{mailto:anthony.roberts@adelaide.edu.au} }}

\begin{document}

\maketitle

\begin{abstract}
Developments in dynamical systems theory provides new support for the discretisation of \pde{}s and other microscale systems.  Here we explore the methodology applied to the gap-tooth scheme in the equation-free approach of Kevrekidis in two spatial dimensions.  The algebraic detail is enormous so we detail computer algebra procedures to handle the enormity.  However,  modelling the dynamics on 2D spatial patches appears to require a mixed numerical and algebraic approach that is detailed in this report.  Being based upon the computation of residuals, the procedures here may be simply adapted to a wide class of reaction-diffusion equations.
\end{abstract}

\tableofcontents

\section{Introduction}

We extend the dynamical systems,  holistic, approach to the macroscale discrete modelling~\cite[e.g.]{MacKenzie00a, Roberts98a, Roberts00a} to two dimensional, homogeneous, nonlinear reaction-diffusion equations on coupled patches of space.    
Following the `equation free' approach of Kevrekidis and colleagues~\cite[e.g.]{Kevrekidis09a}, we address the extraction, using dynamical systems theory, of computationally efficient macroscale models from given microscopic models, whether \pde\ or lattice dynamics or other microscale systems. 
Here we bridge space scales by generalising  to multiple dimensions (specifically~2D) the methodology and supporting theory for the `equation free', gap-tooth method for microsimulators~\cite{Gear03, Samaey03a, Samaey03b}.
As a particular example, this report uses the real valued, two dimensional, Ginzburg--Landau equation 
\begin{equation}
\D tu=\nabla^2u+\alpha(u-u^3)\,. \label{E_gl2d}
\end{equation}
We choose this 2D~real Ginzburg--Landau equation as a prototype \pde\ because it is well studied and its dynamics well understood~\cite[e.g.]{Gibbon93, Levermore96}.  This report details the construction by computer algebra of the macroscale discrete model of its dynamics in two spatial dimensions.  The general theoretical support, the performance and the physical interpretation are detailed elsewhere.

Place the discrete modelling of  two dimensional, reaction-diffusion
equations within the purview of centre manifold theory by dividing the
domain into square patches and introducing special interpatch coupling conditions.  
Define a grid of points~$(X_i,Y_j)$ with, for simplicity, constant spacing~$H$. The $i,j$th patch,~$E_{i,j}$, is then centred upon~$(X_i,Y_j)$ and of width $\Delta x=\Delta y = 2rH$: when the ratio $r<1/2$\,, $E_{i,j}$~forms a patch separated by empty space (gaps) from neighbouring patches; when the ratio $r=1$ the patches overlap and the analysis reduces to that of holistic discretisation~\cite[e.g.]{MacKenzie00a, Roberts98a, Roberts00a}.  Define that $u_{i,j}(x,y,t)$ denotes the field in the $i,j$th~element and so evolves according to the Ginzburg--Landau \pde~\eqref{E_gl2d}.
Using the parameter~$\gamma$ to control the strength of the coupling, use coupling conditions around the $i,j$th~patch of
\begin{equation}
\begin{cases}
u_{i,j}(X_i\pm rH,y,t)= 
\mathcal E_i^{\pm r}(\gamma)\mathcal E_j^{\pm\eta}(\gamma)u_{i,j}(X_i,Y_j,t)\,,
&  |y-Y_j|<rH\,, \\
u_{i,j}(x,Y_j\pm rH,t)= 
\mathcal E_i^{\pm \xi}(\gamma)\mathcal E_j^{\pm r}(\gamma)u_{i,j}(X_i,Y_j,t)\,,
& |x-X_i|<rH\,, 
\end{cases}
\label{eq:epcc}
\end{equation}
In the overlapping element case, $r=1$\,, we use `interelement coupling conditions' around the $i,j$th~element of
\begin{equation}
\begin{cases}
u_{i,j}(X_{i\pm1},y,t)=
\gamma u_{i\pm1,j}(X_{i\pm1},y,t)+(1-\gamma)u_{i,j}(X_i,y,t)
,&  |y-Y_j|<H\,, \\
u_{i,j}(x,Y_{j\pm1},t)=
\gamma u_{i,j\pm1}(x,Y_{j\pm1},t)+(1-\gamma)u_{i,j}(x,Y_j,t),
& |x-X_i|<H\,. 
\end{cases}
\label{EbcsdL}
\end{equation}
These are a natural extension to 2D of those established for 1D dynamics~\cite{Roberts06d}.   Then centre manifold theory~\cite[e.g.]{Carr81, Carr83b, Kuznetsov95} assures us of the existence, relevance and approximation of a slow manifold, macroscale model parametrised by a measure of the amplitude in each element and the coupling strength~$\gamma$.

However,  it appears that models of the \pde~\eqref{E_gl2d} cannot be constructed algebraically.  Thus 
Section~\ref{S_app_2d_num_2ord} details computer algebra to numerically solve for the microscopic subgrid scale field.  
That is, we actually create a coarse grid model of a the fine grid\slash lattice dynamics of an approximation to the \pde.
Its application here to the discretisation of the Ginzburg--Landau \pde~\eqref{E_gl2d} serves as a proof of principle for applying holistic discretisation to general \pde{}s of two or more spatial dimensions.

For example, consider the \pde~\eqref{E_gl2d} when discretised inside patches by a $5\times5$ microscale lattice ($n=2$). Using the patch \ibc~\eqref{eq:epcc} the slow manifold discretisation is
\begin{align}
\dot U_{i,j}={}&
\frac\gamma{H^2}\bdelta^2 U_{i,j}
-\frac{\gamma^2}{12H^2}\big(1-\rat{r^2}{4}\big)\bdelta^4 U_{i,j}
+\frac{\gamma^3}{90H^2}\left(1-\rat{r^2}{4}\right)\left(1-\rat{r^2}{16}\right)\bdelta^6 U_{i,j}
\nonumber\\&{}
+\alpha( U_{i,j}- U_{i,j}^3)
\nonumber\\&{}
+\alpha\gamma r^2\frac1{18}\left\{
\left[3(\mu_x\delta_x U_{i,j})^2+\rat14(\delta_x^2 U_{i,j})^2\right](2+\delta_x^2)
+(\delta_x^2 U_{i,j})^2
\right.\nonumber\\&\qquad\left.{}
+\left[3(\mu_y\delta_y U_{i,j})^2+\rat14(\delta_y^2 U_{i,j})^2\right](2+\delta_y^2)
+(\delta_y^2 U_{i,j})^2
\right\} U_{i,j}
\nonumber\\&{}
+\Ord{\alpha^2+\gamma^4},
\label{eq:glsm9}
\end{align}
where the bold centred difference operator applies in both spatial dimensions,
\begin{eqnarray}
\bdelta^2 U_{i,j}&=&U_{i+1,j}+U_{i-1,j}+U_{i,j+1}+U_{i,j-1}-4 U_{i,j}\,,
\label{eq:bdelta}\\
\bdelta^4 U_{i,j}&=&U_{i+2,j}+U_{i-2,j}+U_{i,j+2}+U_{i,j-2}
\nonumber\\&&{}
-4(U_{i+1,j}+U_{i-1,j}+U_{i,j+1}+U_{i,j-1}) +12 U_{i,j}\,,
\end{eqnarray}

\section{Construct the algebraic slow manifold}
\label{S_app_2ord_2D}

This section provides computer algebra code to generate \emph{analytic} holistic discretisations of reaction-diffusion \pde{}s in 2D. Different \pde{}s are analysed by changing the nonlinear term.  Higher order models are constructed by changing the order of neglected terms. Unfortunately, analytic construction can only be carried out to low order accuracy: Section~\ref{S_app_2d_num_2ord} describes code for a numerical construction that also applies to patches in the gap-tooth scheme.

All code is written in the computer algebra package \textsc{Reduce}.\footnote{\url{http://www.reduce-algebra.com}}

\subsection{Initialisation}

Set some parameters to improve printing of the results.

\<casm2d\><<<
on div; off allfac; on revpri;
factor gam,h,eps,nu,alf;
>>>

\paragraph{Subgrid variables} The subgrid, intra-element, structures are functions of intra-element microscale variables $\verb|xi|=(x-x_i)/h$ and $\verb|yi|=(y-y_j)/h$.

\<casm2d\><<<
depend xi,x;
let df(xi,x)=>1/h;
depend yi,y;
let df(yi,y)=>1/h;
>>>

\paragraph{Parametrise the slow manifold} Parametrise the slow manifold, macroscale, evolution by the evolving grid values~$\verb|u(i,j)|=U_{i,j}$ such that $dU_{i,j}/dt=\verb|gij|$.

\<casm2d\><<<
operator u;
depend u,t;
let df(u(~k,~m),t)=>sub({i=k,j=m},gij);
>>>

The linear, slow subspace, approximation is that of piecewise constant fields and no evolution.

\<casm2d\><<<
uij:=u(i,j);
gij:=0;
>>>

\paragraph{Set asymptotic truncation}  Here scale nonlinearity parameter~$\alpha$ with parameter~$\gamma$ so we truncate to residuals and errors of order~$\Ord{\alpha^3+\gamma^3}$.  

\<casm2d\><<<
let eps^3=>0; 
gamma:=gam*eps;
alpha:=alf*eps;
>>>

\paragraph{General multinomial} For the method of undetermined coefficients we set up a general multinomial solution with unknown coefficients~\verb|cc(m,n)|, up to order~\verb|o| in the intra-element variables.  Collect the unknown coefficients in the set~\verb|cs|.

\<casm2d\><<<
o:=6;
operator cc;
vv:=for m:=0:o sum for n:=0:o-m sum cc(m,n)*xi^m*yi^n$
cs:=for m:=0:o join for n:=0:o-m collect cc(m,n)$
>>>

Need to find the coefficients in multivariate expressions so define a recursive procedure \verb|coeffs| that given a list of expressions and a list of variables it creates a new extended list of all the coefficients of the variables.

\<casm2d\><<<
procedure coeffs(exps,vars); 
    if vars={} then exps else coeffs(
        foreach z in exps join coeff(z,first vars),
        rest vars);
>>>

\subsection{Iteratively construct the slow manifold}

Now start the iteration, repeating corrections until all residuals are zero (to a maximum of ten iterations).

\<casm2d\><<<
for it:=1:10 do begin 
>>>

\paragraph{Compute residuals of governing equations}
Compute the residual of the reaction-diffusion \pde\ in the general $i,j$th~element.  Could easily change the reaction term to any polynomial in the field~$u_{i,j}$.  Could also incorporate advection terms into the microscale \pde.
Also compute the residuals of the inter-element coupling \ibc{}s.

\<eqnResiduals\><<<
de:=df(uij,t)-df(uij,x,2)-df(uij,y,2)-alpha*(uij-uij^3)$
bcr:=sub(xi=1,uij)-sub(xi=0,uij)
   -gamma*(sub({xi=0,i=i+1},uij)-sub(xi=0,uij))$
bcl:=sub(xi=-1,uij)-sub(xi=0,uij)
   -gamma*(sub({xi=0,i=i-1},uij)-sub(xi=0,uij))$
bct:=sub(yi=1,uij)-sub(yi=0,uij)
   -gamma*(sub({yi=0,j=j+1},uij)-sub(yi=0,uij))$
bcb:=sub(yi=-1,uij)-sub(yi=0,uij)
   -gamma*(sub({yi=0,j=j-1},uij)-sub(yi=0,uij))$
>>>

Include the computation of the residuals both here and later.

\<casm2d\><<<
`<eqnResiduals`>
>>>

For information as to the progress of the iteration, print out the lengths of the residuals: when all are one, then all residuals are probably zero.

\<casm2d\><<<
write lengths:={length(de),length(bcr),length(bcl)
    ,length(bct),length(bcb)};
>>>

\paragraph{Add the as yet unknown corrections}
To find the desired update in each iteration, first substitute the form of the update into the computed residuals.

\<casm2d\><<<
deq:=de+gd-df(vv,x,2)-df(vv,y,2);
rbcr:=bcr+sub(xi=1,vv)-sub(xi=0,vv);
lbcl:=bcl+sub(xi=-1,vv)-sub(xi=0,vv);
tbct:=bct+sub(yi=1,vv)-sub(yi=0,vv);
bbcb:=bcb+sub(yi=-1,vv)-sub(yi=0,vv);
>>>

\paragraph{Solve for the corrections}
Then extract equate coefficients of each power of the multinomial in the intra-element variables, and append boundary conditions.

\<casm2d\><<<
eqns:=append(coeffs({deq},{xi,yi}),cc(0,0). 
   append(coeff(rbcr,yi),append(coeff(lbcl,yi),
   append(coeff(tbct,xi),coeff(bbcb,xi)))));
sol:=solve(eqns,gd.cs);
>>>

\paragraph{Update}
Update the field and the evolution, assuming a solution was found (not true for higher orders).

\<casm2d\><<<
uij:=uij+sub(sol,vv);
gij:=gij+sub(sol,gd);
>>>

\paragraph{Terminate}
End the iteration when all residuals are zero, or too many iterations have been performed.  

\<casm2d\><<<
showtime;
if {de,bcr,bcl,bct,bcb}={0,0,0,0,0} then write it:=1000000+it;
end;
>>>

\subsection{Scrounge an extra order of evolution using solvability}

First define the linear integral operator $\verb|inthat(a,xi)|=\int_{-1}^1 (1-|\xi|)\verb|a|\,d\xi$ in order to quickly apply the solvability condition.

\<casm2d\><<<
operator inthat;
linear inthat;
let { inthat(~a^~p,~b)=>0 when (a=b)and(not evenp(p))
, inthat(~a^~p,~b)=>2/(p+1)/(p+2) when (a=b)and evenp(p)
, inthat(~a,~b)=>0 when (a=b)
, inthat(1,~b)=>1 
};
>>>

Set the requisite \emph{next order} of truncation in the asymptotic expansion by finding the highest order currently retained in coupling~$\gamma$, then setting to discard terms of two orders higher.  Use the `instant evaluation' property of for-loop variables in Reduce.

\<casm2d\><<<
o:=2+deg((1+eps)^9,eps)$
for p:=o:o do let eps^p=>0;
>>>

Compute exactly the same residuals of the \pde\ and the inter-element coupling conditions.

\<casm2d\><<<
`<eqnResiduals`>
>>>

Compute the next correction to the evolution~\verb|gij| by integrating the residual of the \pde\ over an element, and including the contributions from the boundary coupling residuals.

\<casm2d\><<<
gd:=inthat(de,xi)$
gd:=inthat(gd,yi)$
gd:=gd+inthat(bcr+bcl,yi)/h^2+inthat(bct+bcb,xi)/h^2;
gij:=gij-gd$
showtime;
>>>

\section{Postprocessing analytic}

\subsection{Write a sci/matlab function script}

Optionally generate an efficient sci/matlab function in file~\verb|gl2dnumred|.  This function file is to be processed with the unix script~\verb|red2mat| and then it can be used in matlab to integrate solutions using \ode\ solvers. 

\<casm2d\><<<
write "
Generating mat/scilab function file in gl2dnumred
";
`<functionfile`>
write "finished generating mat/scilab file";
write "Use sed -f red2mat gl2dnumred > gl2dnpde.m";
showtime;
>>>

Improve printing.
\<functionfile\><<<
linelength 66$ 
off nat$
>>>

Need boundary conditions that depend upon the width of the stencil of the discretisation.  Get this width from the truncation in parameter~\verb|eps|, assuming coupling~$\gamma$ is proportional to~\verb|eps|.  If not, then have to change. \<functionfile\><<<
o:=deg((1+eps)^9,eps);
>>>

Want code for the fully coupled system, $\gamma=1$, and without the artificial parameter~\verb|eps|.
\<functionfile\><<<
gam:=eps:=1;
>>>

Use the scope package to generate efficient computation of the discretisation.
\<functionfile\><<<
load_package scope;
>>>

Now direct output to the file \verb|gl2dnumred|.
\<functionfile\><<<
out "gl2dnumred"$
>>>

Write the preamble for the function.  Given a vector of macroscale variables, unpack them into an array of 2D grid values.  Get the grid size and the gap-tooth ratio as global variables.

\<functionfile\><<<
write "function udot=gl2dnpde(t,uv);
global h rr;
nn=length(uv);
n=sqrt(nn);
u=zeros(n,n); u(:)=uv(1:nn);"$
>>>

Allow for two possible boundary conditions: either doubly periodic; or odd-even solutions.  Choose in the reduce code by changing the switch in this if-statement.

\<functionfile\><<<
if 0 then 
`<periodicBCs`>else  `<oddevenBCs`>$
>>>

Pad the data array with the periodic BCs: the two cases depend upon the stencil width.

\<periodicBCs\><<<
if o=1 then write "
u=[u(end,:);u;u(1,:)];
u=[u(:,end),u,u(:,1)]"
else write "
u=[u(end-",o-1,":end,:);u;u(1:",o,",:)];
u=[u(:,end-",o-1,":end),u,u(:,1:",o,")]"
>>>

Alternatively, pad the data array with the odd/even BCs: the two cases depend upon the stencil width.  The odd/even BCs suit a sine expansion in~$x$ and a cosine expansion in~$y$.

\<oddevenBCs\><<<
if o=1 then write "
u=[-u(1,:);u;-u(end,:)];
u=[+u(:,1),u,+u(:,end)]"
else write "
u=[-u(",o,":-1:1,:);u;-u(end:-1:end-",o-1,",:)];
u=[+u(:,",o,":-1:1),u,+u(:,end:-1:end-",o-1,")]"
>>>

Define index vectors that apply to every macroscale grid point in the interior of the padded data array.
\<functionfile\><<<
write "j=",o,"+(1:n); i=j"$
>>>

Use the scope package optimisation to generate the code for the computation of an array of time derivatives.
\<functionfile\><<<
optimize udot:=:gij iname o$
>>>

Write code to form the array into a vector.
\<functionfile\><<<
write "udot=udot(:)"$
>>>

Close the function file as all code is done.
\<functionfile\><<<
shut "gl2dnumred";
>>>

Undo a printing option, and clear \verb|gam| and \verb|eps| so we can see the coupling parameter influence again.
\<functionfile\><<<
on nat;
clear gam;
clear eps;
>>>

\subsection{Derive the equivalent PDE}

Optionally derive the equivalent \pde\ for this model.

\<casm2d\><<<
if 1 then begin write "
Optionally find the equivalent partial differential equation
of the holistic discretisation.
";
`<equivpde`>
    showtime;
end;
>>>

Depending upon the order of inter-element coupling, as measured by~\verb|eps|, then set the order of truncation in the grid size~\verb|h|.  This code assumes \verb|gamma| is first order in~\verb|eps|, and needs changing if otherwise.  Use a dummy for-loop because, in Reduce. loop parameters, here~\verb|q|, have the ``instant evaluation'' property that means it can be used usefully as an index in the \verb|let| statement.

\<equivpde\><<<
    o:=2+2*deg((1+eps)^9,eps);
    for q:=o:o do let h^(q+1)=>0;
>>>

Use operator \verb|uu(p,q)| to denote the mixed derivative $\partial^{p+q}U/\partial x^p\partial y^q$.  Hence define a rule to expand \verb|u(i+k',j+l')| in a multivariable Taylor series.  Reduce does not like~\verb|0^0| so handle the zero cases separately in the summation.

\<equivpde\><<<
    operator uu;
    rules:={u(~k,~l)=> uu(0,0)
        +(for q:=1:o sum uu(0,q)*(l-j)^q*h^q/factorial(q))
        +(for p:=1:o sum uu(p,0)*(k-i)^p*h^p/factorial(p))
        +(for p:=1:o sum for q:=1:o-p sum
          uu(p,q)*(k-i)^p*(l-j)^q*h^(p+q)/factorial(p)/factorial(q))
        }$
>>>

Make different printing.
\<equivpde\><<<
    off revpri; on list;
>>>

Transform the discrete model to its equivalent \pde using the Taylor expansion.  Remove the artificial expansion parameter~\verb|eps| by settng to one. Sum at full coupling~\verb|gam|; remove this to see the~$\gamma$ dependence.

\<equivpde\><<<
    write
    gpde:=(sub({eps=1,gam=1},gij) where rules)$
    off list;
>>>

\subsection{Derive the centred finite difference form}

Optionally derive the finite difference form of the evolution.

\<casm2d\><<<
if 1 then begin write "
Generate evolution in terms of centred finite difference
operators md=mu*delta and dd=delta^2 in each of the x and y
directions.
";
`<finitediff`>
    showtime;
end;
>>>

Define the four centred difference operators, two in each spatial direction.  Then define commutation rules so that expressions involving these operators change to a canonical form where $x$~operators come before $y$~operators, and $\mu\delta$~operators comes before $\delta^2$~operators.  Lastly, invoke the operator rule that $\mu^2=1+\delta^2/4$.

\<finitediff\><<<
operator mdx; operator mdy;
operator ddx; operator ddy;
linear mdx; linear mdy;
linear ddx; linear ddy;
let { mdy(mdx(~a,~u),~v)=>mdx(mdy(a,v),u)
    , ddx(mdx(~a,~u),~v)=>mdx(ddx(a,v),u)
    , ddy(mdx(~a,~u),~v)=>mdx(ddy(a,v),u)
    , ddx(mdy(~a,~u),~v)=>mdy(ddx(a,v),u)
    , ddy(mdy(~a,~u),~v)=>mdy(ddy(a,v),u)
    , ddy(ddx(~a,~u),~v)=>ddx(ddy(a,v),u)
    , mdx(mdx(~a,~u),~u)=>ddx(a,u)+ddx(ddx(a,u),u)/4
    , mdy(mdy(~a,~u),~u)=>ddy(a,u)+ddy(ddy(a,u),u)/4
    };
>>>

Define the shift rules that $U_{i\pm1}=U_i\pm\mu\delta U_i+\rat12\delta^2U_i$. Apply them to the four cases of $U_{i\pm k,j}$ and~$U_{i,j\pm k}$.

\<finitediff\><<<
rules:={u(i,j)=>u
    ,u(i+~k,~j)=>u(i+k-1,j)+mdx(u(i+k-1,j),u)
        +1/2*ddx(u(i+k-1,j),u) when k>0
    ,u(i+~k,~j)=>u(i+k+1,j)-mdx(u(i+k+1,j),u)
        +1/2*ddx(u(i+k+1,j),u) when k<0
    ,u(~i,j+~k)=>u(i,j+k-1)+mdy(u(i,j+k-1),u)
        +1/2*ddy(u(i,j+k-1),u) when k>0
    ,u(~i,j+~k)=>u(i,j+k+1)-mdy(u(i,j+k+1),u)
        +1/2*ddy(u(i,j+k+1),u) when k<0
};
>>>

Transform the evolution of the grid values.  Could also transform the microscale structure in~$\uij$, but have not done so here.

\<finitediff\><<<
write gop:=(sub(eps=1,gij) where rules);
>>>

\subsection{Fin analytic}

Give the overall grand final `end' statement, then output all the code.

\<casm2d\><<<
end;
>>>

\OutputCode\<casm2d\>

\section{Numerically construct the slow manifold}
\label{S_app_2d_num_2ord}

This section describes computer algebra code to generates macroscale models of the gap-tooth scheme applied to reaction-diffusion \pde{}s in 2D. 
Different \pde{}s are analysed by simply changing the nonlinear term, or more carefully changing the Laplacian.  
Higher order models are constructed by changing the order of neglected terms.  
Instead of trying to solve analytically for the subgrid scale structure, here solve for the structure numerically.

All code is written in the freely available computer algebra package \textsc{Reduce}.\footnote{\url{http://www.reduce-algebra.com}}

\subsection{Initialisation}
Improve the format of printing.

\<ncsm2d\><<<
on div; off allfac; on revpri;
factor gam,h,eps,nu,alf;
>>>

Load routines to do the $LU$~decomposition, Section~\ref{sec:lu}, and subsequent back substitutions, Section~\ref{sec:bs}.

\<ncsm2d\><<<
in "lu_decomp.red"$ 
in "lu_backsub.red"$
>>>

Need to find the coefficients in multivariate expressions so define a recursive procedure \verb|coeffs| that given a list of expressions and a list of variables it creates a new extended list of all the coefficients of the variables.

\<ncsm2d\><<<
procedure coeffs(exps,vars); 
    if vars={} then exps else coeffs(
        foreach z in exps join coeff(z,first vars),
        rest vars);
>>>

Use patches of size ratio~$r$ where: $r=1$ is holistic overlapping; $r=1/2$ the elements abut; and for $r<1/2$ there is a gap in the gap-tooth scheme.

\<ncsm2d\><<<
>>>

The subgrid, microscale, numerical resolution is improved by increasing~\verb|n|.  The subgrid grid step is $\verb|dx|=H/n$ and the number of equations and unknowns is $\verb|nn|=N=(2n+1)^2+1$\,.  Takes some minutes to construct solutions for lattice $n=8$, but it is reasonably doable; have not tried $n\geq9$ yet.

\<ncsm2d\><<<
n:=2; 
dx:=r*h/n$ 
ns:=2*n+1$ 
nn:=ns^2+1$ 
>>>

Define the matrices used to store the subgrid, and to represent and solve the equations.  The scope of matrices are global in Reduce.

\<ncsm2d\><<<
matrix eqns(nn,1),zeros(nn,1);
matrix indx(nn,1),vv(nn,1),lu(nn,nn); 
>>>

\paragraph{Parametrise the slow manifold} The slow manifold of the macroscale discretisation is to be parametrised by the evolution of $u_{i,j}=\verb|u(i,j)|$.   The evolution $dU_{i,j}/dt=g_{i,j}=\verb|gij|$ and unknown updates are stored in \verb|gd|.

\<ncsm2d\><<<
operator u; depend u,t;
let df(u(~k,~m),t)=>sub({i=k,j=m},gij);
gij:=gd/dx^2; 
>>>

Initialise with constant field in $\verb|uij|=u_{i,j}$.  Also add in the unknown update field~\verb|ud|.  Use ii and jj as the subgrid microscale lattice indices.

\<ncsm2d\><<<
matrix uij(ns,ns)$
operator ud;
for ii:=1:ns do for jj:=1:ns do 
    uij(ii,jj):=u(i,j)+ud(ii,jj)$
>>>

\paragraph{Possibly include a forcing function}  \verb|fu(p,q)=df(ff,x,p,y,q)|, the physical derivatives, evaluated at the $(i,j)$ grid point.  But omit for the moment.

\<ncsm2d\><<<
operator fu;
maxf:=0;
ffij:=(for m:=0:maxf sum h^m
    *(for n:=0:m sum fu(n,m-n)*xi^n*yi^(m-n)));
>>>

\paragraph{Set asymptotic truncation}  Truncate the asymptotic expansion in powers of the coupling parameter~$\gamma=\verb|gam|$.  Here scale the nonlinearity parameter $\alpha$ with~$\gamma$ to most easily control the truncation.  

\<ncsm2d\><<<
let { eps^4=>0 }; 
gamma:=gam*eps;
gammad:=1-gamma$
alpha:=alf*eps^2; 
nu:=1; 
>>>

Also ignore the `updates' when multiplied by the small parameter so that the equations remain linear in the updates.

\<ncsm2d\><<<
let { eps*ud(~i,~j)=>0, eps*gd=>0,
      gam*ud(~i,~j)=>0, gam*gd=>0 }; 
>>>

\subsection{Classic Lagrange interpolation}

Coupling conditions are written in terms of classic Lagrange interpolation.  Thus set up the interpolation fields here as functions of the surrounding macroscale grid values.  As written so far we can go to errors~\Ord{\gamma^5}.  Suspect I should be able to code this better.

First define procedures for centred mean and difference operators in each of the two spatial directions.

\<ncsm2d\><<<
procedure mudi(u); (sub(i=i+1,u)-sub(i=i-1,u))/2;
procedure mudj(u); (sub(j=j+1,u)-sub(j=j-1,u))/2;
procedure dedi(u); (sub(i=i+1,u)-2*u+sub(i=i-1,u));
procedure dedj(u); (sub(j=j+1,u)-2*u+sub(j=j-1,u));
>>>

Second, interpolate in the $x$~direction where $\verb|xi|=(x-X_{i,j})/H$.

\<ncsm2d\><<<
uu:=u(i,j)$
uu:=uu+gamma*(+mudi(uu)+dedi(uu)*xi/2)*xi
      +gamma^2*dedi(+mudi(uu)+dedi(uu)*xi/4)*xi*(xi^2-1)/6
      +gamma^3*dedi(dedi(+mudi(uu)+dedi(uu)*xi/6))
      *xi*(xi^2-1)*(xi^2-4)/120
      +gamma^4*dedi(dedi(dedi(+mudi(uu)+dedi(uu)*xi/8)))
      *xi*(xi^2-1)*(xi^2-4)*(xi^2-9)/5040
      $
>>>

Third, interpolate in the $y$~direction where $\verb|yi|=(y-Y_{i,j})/H$.

\<ncsm2d\><<<
uu:=uu+gamma*(+mudj(uu)+dedj(uu)*yi/2)*yi
      +gamma^2*dedj(+mudj(uu)+dedj(uu)*yi/4)*yi*(yi^2-1)/6
      +gamma^3*dedj(dedj(+mudj(uu)+dedj(uu)*yi/6))
      *yi*(yi^2-1)*(yi^2-4)/120
      +gamma^4*dedj(dedj(dedj(+mudj(uu)+dedj(uu)*yi/8)))
      *yi*(yi^2-1)*(yi^2-4)*(yi^2-9)/5040
      $
>>>

\subsection{Iteratively construct the slow manifold}

Repetitively find updates to the subgrid microscale field until all residuals are zero.

\<ncsm2d\><<<
for iter:=1:19 do begin
>>>

\paragraph{Update the u field and evolution} 
Form the residuals of lattice equations and all the coupling conditions, using~\verb|q| to count the number of equations.
First the corner equations.

\<ncsm2d\><<<
    write "Compute the u residual";
    q:=0; 
    eqns(q:=q+1,1):=ud(1,1);
    eqns(q:=q+1,1):=ud(1,ns);
    eqns(q:=q+1,1):=ud(ns,1);
    eqns(q:=q+1,1):=ud(ns,ns);
>>>

Now adjoin the residuals of the interpatch\slash element coupling conditions.  First, if $r=1$, then use the interelement coupling conditions.

\<ncsm2d\><<<
    if r=1 then for ll:=2:ns-1 do begin
      eqns(q:=q+1,1):=uij(ns,ll)-(1-gamma)*uij(n+1,ll)
          -gamma*sub(i=i+1,uij(n+1,ll)); 
      eqns(q:=q+1,1):=uij(1,ll)-(1-gamma)*uij(n+1,ll)
          -gamma*sub(i=i-1,uij(n+1,ll)); 
      eqns(q:=q+1,1):=uij(ll,ns)-(1-gamma)*uij(ll,n+1)
          -gamma*sub(j=j+1,uij(ll,n+1)); 
      eqns(q:=q+1,1):=uij(ll,1)-(1-gamma)*uij(ll,n+1)
          -gamma*sub(j=j-1,uij(ll,n+1)); 
    end
>>>

For the following patch coupling conditions (when $r<1$), recall that \verb|uu| stores an expression for classic Lagrange interpolation from neighbouring grid values.

\<ncsm2d\><<<
    else for ll:=2:ns-1 do begin
      xy:=r*(ll-n-1)/n; 
      eqns(q:=q+1,1):=uij(ns,ll)-sub({xi=r,yi=xy},uu); 
      eqns(q:=q+1,1):=uij(1,ll)-sub({xi=-r,yi=xy},uu); 
      eqns(q:=q+1,1):=uij(ll,ns)-sub({yi=r,xi=xy},uu); 
      eqns(q:=q+1,1):=uij(ll,1)-sub({yi=-r,xi=xy},uu); 
    end;
>>>

Adjoin the residuals of the interior subgrid equations which are a simple, second order, centred difference approximations to the \pde. 

\<ncsm2d\><<<
    for ii:=2:ns-1 do for jj:=2:ns-1 do
    eqns(q:=q+1,1):=-df(uij(ii,jj),t)*dx^2
        +nu*((uij(ii+1,jj)-2*uij(ii,jj)+uij(ii-1,jj))
            +(uij(ii,jj+1)-2*uij(ii,jj)+uij(ii,jj-1)))
        +alpha*dx^2*(uij(ii,jj)-uij(ii,jj)^3);
>>>

Last is the amplitude condition that the central value is the unchanged amplitude (could make the amplitude equal to the element mean if necessary), and set \verb|ok| if all the right-hand sides are zero.

\<ncsm2d\><<<
    eqns(q:=q+1,1):=ud(n+1,n+1);
    ok:=if eqns=zeros then 1 else 0; 
>>>

\paragraph{Form matrix of the linear equations}
In the first iteration only, now form matrix of the linear equations from the perturbation variables \verb|ud| and \verb|gd| that have been carried around in this first iterate. Extract these unknowns and assign their coefficients to the \verb|lu|~matrix.  Then perform the $LU$~decomposition for later solution of the homological equation.  Lastly remove the symbolic unknowns.

\<ncsm2d\><<<
    if iter=1 then begin
        for qq:=1:nn do begin q:=0;
            for ii:=1:ns do for jj:=1:ns do
            lu(qq,q:=q+1):=-coeffn(eqns(qq,1),ud(ii,jj),1);
            lu(qq,q:=q+1):=-coeffn(eqns(qq,1),gd,1);
            end;
        write "starting LU decomposition to store";
        lu_decomp();
        write "finished LU decomposition";
        showtime;
        Comment get rid of perturbation variables hereafter;
        let { ud(~ii,~jj)=>0, gd=>0 };
    end; 
>>>

Use the factorisation stored in~\verb|lu| to solve for the updates, and assign to the relevant $g$~and~$u_{i,j}$.

\<ncsm2d\><<<
    write "Back substitution solves for updates";
    lu_backsub();
    gij:=gij+eqns(nn,1)/dx^2;
    q:=0;
    for ii:=1:ns do for jj:=1:ns do
        uij(ii,jj):=uij(ii,jj)+eqns(q:=q+1,1);
>>>

Terminate the loop when all residuals are zero.  

\<ncsm2d\><<<
    showtime;
    if ok then write iter:=1000000+iter;
end;
>>>

\section{Postprocessing numerical}

\subsection{Write a sci/matlab function script}

Optionally generate an efficient sci/matlab function in file~\verb|gl2dnumred|.  This function file is to be processed with the unix script~\verb|red2mat| and then it can be used in matlab to integrate solutions using \ode\ solvers. 

\<ncsm2d\><<<
write "
Generating mat/scilab function file in gl2dnumred
";
`<functionfile`>
write "finished generating mat/scilab file";
write "Use sed -f red2mat gl2dnumred > gl2dnpde.m";
showtime;
>>>

\subsection{Derive the equivalent PDE}

Optionally derive the equivalent \pde\ for this model.

\<ncsm2d\><<<
if 1 then begin write "
Optionally find the equivalent partial differential equation
of the holistic discretisation.
";
`<equivpde`>
    showtime;
end;
>>>

\subsection{Derive the centred difference form}

Optionally derive the finite difference form of the evolution.

\<ncsm2d\><<<
if 1 then begin write "
Generate evolution in terms of centred finite difference
operators md=mu*delta and dd=delta^2 in each of the x and y
directions.
";
`<finitediff`>
    showtime;
end;
>>>

Additionally write the coefficients of three parts in the $\Ord{\alpha^3+\gamma^3}$ model.

\<finitediff\><<<
write g2nd:=factorize(coeffn(h^2*gop,gam,2));
write g3rd:=factorize(coeffn(h^2*gop,gam,3));
write gnon:=factorize(coeffn(coeffn(h^2*gop,alf,1),gam,1));
>>>

Gather results of \verb|gnon| into a list and fit rational Pade functions in~$1/n^2$.  Need an odd number of results.  Note that $\verb|rn|=1/n$ and the multiplication by four of the coefficient.

\<finitediff\><<<
write "Some meta-analysis of different subgrid resolutions";
ns:={2,3,4,5,6,7,8};
cs:={1/72
    ,1/60
    ,179/10136
    ,775/42762
    ,679909/36998632
    ,237808723/12834019900
    ,133046058951/7141880630840
    };
operator a; operator b;
procedure paden(ns,cs);
begin
    o:=(length(ns)-1)/2;
    vars:=for p:=0:o join {a(p),b(p)};
    eqns:=(b(0)=1).(for l:=1:length(ns) collect 
        (for p:=0:o sum b(p)/part(ns,l)^(2*p))*part(cs,l)*4=
        (for p:=0:o sum a(p)/part(ns,l)^(2*p)) );
    soln:=solve(eqns,vars);
    return sub(soln,(for p:=0:o sum a(p)*rn^(2*p))
                            /(for p:=0:o sum b(p)*rn^(2*p)));
end;
padefull:=paden(ns,cs);
on rounded; write largenlimit:=sub(rn=0,padefull); off rounded;
write pade4:=paden({2,3,4},
    {1/72
    ,1/60
    ,179/10136});
>>>

\subsection{Fin numerical}

Give the grand final `end' statement, and output all the code.

\<ncsm2d\><<<
end;
>>>

\OutputCode\<ncsm2d\>

\section{Factorisation and solution of linear equations}

\subsection{$LU$ decomposition}
\label{sec:lu}

This $LU$~decomposition was adapted from Numerical recipes in Fortran~77~\cite{Press92}.  Requires matrices \verb|lu(n,n)|, \verb|indx(n,1)| and \verb|eqns(n,1)| to be predefined. Matrices are global in scope in Reduce so we do not bother passing parameters.

\<lu_decomp\><<<
procedure lu_decomp;
begin scalar n,np,i,j,k,imax,dd,aamax,dum,sum,tiny;
  tiny:=1.0e-20;
  n:=first(length(lu));
  dd:=1;
  for i:=1:n do begin
    aamax:=abs(lu(i,1));
    for j:=2:n do aamax:=max(aamax,abs(lu(i,j)));
    vv(i,1):=1/aamax;
  end;
  for j:=1:n do begin
    for i:=1:(j-1) do 
      lu(i,j):=lu(i,j)-(for k:=1:(i-1) sum lu(i,k)*lu(k,j));
    aamax:=0;
    for i:=j:n do begin
      lu(i,j):=lu(i,j)-(for k:=1:(j-1) sum lu(i,k)*lu(k,j));
      dum:=vv(i,1)*abs(lu(i,j));
      if dum>=aamax then begin
        imax:=i;
        aamax:=dum;
      end;
    end;
    if not(j=imax) then begin
      for k:=1:n do begin
        dum:=lu(imax,k);
        lu(imax,k):=lu(j,k);
        lu(j,k):=dum;
      end;
      dd:=-dd;
      vv(imax,1):=vv(j,1);
    end;
    indx(j,1):=imax;
    if lu(j,j)=0 then lu(j,j):=tiny;
    if not(j=n) then begin
      dum:=1/lu(j,j);
      for i:=(j+1):n do lu(i,j):=lu(i,j)*dum;
    end;
  end;
end;
end;
>>>

\OutputCode\<lu_decomp\>

\subsection{$LU$ back substitution}
\label{sec:bs}

This $LU$~back substitution was adapted from Numerical recipes in Fortran~77~\cite{Press92}.  Requires matrices \verb|lu(n,n)|, \verb|indx(n,1)| and \verb|eqns(n,1)| to be predefined: the first two must be obtained from procedure \verb|lu_decomp|; and the last may be algebraic. Matrices are global in scope in Reduce so we do not bother passing parameters.

\<lu_backsub\><<<
procedure lu_backsub;
begin scalar n,i,j,ii,ll,sum;
  n:=first(length(lu));
  ii:=0;
  for i:=1:n do begin
    ll:=indx(i,1);
    sm:=eqns(ll,1);
    eqns(ll,1):=eqns(i,1);
    if not(ii=0) then 
       sm:=sm-(for j:=ii:(i-1) sum lu(i,j)*eqns(j,1))
    else if not(sm=0) then ii:=i;
    eqns(i,1):=sm;
  end;
  for i:=n step -1 until 1 do eqns(i,1):=(eqns(i,1)
     -(for j:=(i+1):n sum lu(i,j)*eqns(j,1)))/lu(i,i);
end;
end;
>>>

\OutputCode\<lu_backsub\>

\paragraph{Acknowledgement}  The Australian Research Council Discovery Project grants DP0774311 and DP0988738 helped support this research.

\bibliographystyle{plain}
\bibliography{more2,ajr,bib}

\end{document}